\documentclass[11pt,reqno]{amsart}
\usepackage{amsmath,amsfonts, amssymb, amsthm}
  \usepackage{graphics} 
  \usepackage{epsfig} 
\usepackage{graphicx}  \usepackage{epstopdf}
 \usepackage[colorlinks=true]{hyperref}
\hypersetup{urlcolor=blue, citecolor=red}

\newtheorem{conjecture}{Conjecture}[section]

\theoremstyle{definition}

\newtheorem{remark}{Remark}[section]

\def\pmod #1{\ ({\rm{mod}}\ #1)}
\def\Z{\Bbb Z}
\def\N{\Bbb N}

\def\l{\left}
\def\r{\right}
\def\bg{\bigg}
\def\({\bg(}
\def\){\bg)}
\def\t{\text}
\def\f{\frac}

\def\mo{{\rm{mod}\ }}
\def\pmod#1{\ (\mo\ #1)}

\def\ls{\leqslant}
\def\gs{\geqslant}

\def\sm{\setminus}
\def\bi{\binom}

\def\eq{\equiv}

\begin{document}
\hbox{Preprint, {\tt arXiv:2603.29973}}
\medskip

\title[Various conjectural series identities]
      {Various conjectural series identities}
\author[Zhi-Wei Sun]{Zhi-Wei Sun}


\address{School of Mathematics, Nanjing
University, Nanjing 210093, People's Republic of China}
\email{{\tt zwsun@nju.edu.cn}
\newline\indent
{\it Homepage}: {\tt http://maths.nju.edu.cn/\lower0.5ex\hbox{\~{}}zwsun}}

\keywords{Binomial coefficients, combinatorial identities, infinite series, Riemann's zeta function, congruences.
\newline \indent Supported by the Natural Science Foundation of China (grant no. 12371004).}

\begin{abstract} In this paper we collect over 150 new series identities
(involving binomial coefficients) conjectured by the author in  2026. The values involved
 are related to $\pi$ or Riemann's zeta function or Dirichlet's $L$-function.
For example, we conjecture that
$$\sum_{k=0}^\infty\frac{16k+3}{(-202^2)^k}
 \binom{2k}kT_k(19,-20)T_{2k}(9,-5)=\frac{43\sqrt{101}}{75\pi},$$
 where $T_n(b,c)$ denotes the coefficient of $x^n$ in the expansion of $(x^2+bx+c)^n$.
 The conjectures in this paper might interest some readers and stimulate further research.
\end{abstract}
\maketitle

\section{Introduction}
\setcounter{equation}{0}
 \setcounter{conjecture}{0}
 \setcounter{theorem}{0}
 \setcounter{proposition}{0}


 The classical Ramanujan-type series for $1/\pi$ (cf. Cooper \cite[Chapter 14]{Co17} and Ramanujan \cite{R}) takes the form:
 \begin{equation}\label{R}\sum_{k=0}^\infty(ak+b)\f{\bi{2k}kc_k}{m^k}=\f{\sqrt d}{\pi},\end{equation}
 where $a,b,m$ are integers with $am\not=0$, $d\in\Z^+=\{1,2,3,\ldots\}$ and $c_k$ takes one of the forms
 $$\bi{2k}k^2,\ \bi{2k}k\bi{3k}k,\ \bi{2k}k\bi{4k}{2k},\ \bi{3k}k\bi{6k}{3k}.$$

 For $b,c\in\Z$ and $n\in\N$, the generalized central trinomial coefficient
 $T_n(b,c)$ is the coefficient of $x^n$ in the expansion of $(x^2+bx+c)^n$.
 As $T_n(2,1)=\bi{2n}n$, we may view $T_n(b,c)$ as a natural extension of the central binomial coefficient $\bi{2n}n$. Motivated by this and Ramanujan-type series, the author posed conjectural series for $1/\pi$
 involving generalized central trinomial coefficients of types I-VII in 2011
 (cf. \cite{S14b}), type VIII in 2019 (cf. \cite{ERA}), and type IX in 2020 (cf. \cite{Colloq}).

 In this paper we collect the author's various open conjectural series identities found in 2026.

In the next section, we pose three series for $1/\pi$ of type X involving generalized central trinomial coefficients.

Recall that the harmonic numbers are those rational numbers
 $$H_n:=\sum_{0<k\ls n}\f1k\ \ (n=0,1,2,\ldots).$$
 For any integer $m>1$, the harmonic numbers of order $m$ are given by
 $$H_n^{(m)}:=\sum_{0<k\ls n}\f1{k^m}\ \ (n=0,1,2,\ldots).$$
 The author \cite[Conjecture 29]{harmonic} guessed that each Ramanujan-type series for $1/\pi$ has an explicit variant with summands involving certain linear combination of harmonic numbers, and the conjecture was confirmed by Yajun Zhou \cite{ZhouRJ}.
 For many open conjectures on series involving (generalized) harmonic numbers
 discovered before 2026, see a series of the author's papers
 \cite{Sab,S24,harmonic,III,II}. See also \cite{GL,HHW,SunZhou,Wei,WR} for some results on series involving harmonic numbers.

 In Section 3, we pose over ten new conjectures on series involving harmonic numbers.

The Riemann zeta function is given by
$$\zeta(s)=\sum_{k=1}^\infty \f1{k^s}\ \ \t{for}\ \Re(s)>1.$$
If we set $f(x)=x^{-s}=e^{-s\log x}$ for $x>0$, then for any positive integer $m$ we have
$$\sum_{k=1}^\infty f^{(m)}(k)=\sum_{k=1}^\infty \f{\prod_{j=0}^{m-1}(-s-j)}{k^{m+s}}
=(-1)^m\zeta(s+m)\prod_{j=0}^{m-1}(s+j),$$
where $f^{(m)}(x)$ denotes the $m$th order derivative of $f(x)$.

Recall that the Gamma function is given by
$$\Gamma(z)=\int_0^\infty t^{z-1}e^{-t}dt\  \t{for}\ \Re(z)>0.$$
For integers $m>n>0$, the binomial coefficient $\bi {mk}{nk}$ with $k\in\N=\{0,1,2,\ldots\}$
is just the value of the function
$$\f{\Gamma(mx+1)}{\Gamma(nx+1)\Gamma((m-n)x+1)}$$
at $x=k$, thus $\bi{mk}{nk}$ is meaningful whenever $k>-1/m$.
For a Ramanujan-type series $\sum_{k=0}^\infty f(k)$ of the form \eqref{R}, we know the value of 
$\sum_{k=0}^\infty f'(k)$ by Sun \cite[Conjecture 29]{harmonic} as confirmed by Zhou \cite{ZhouRJ}.
For example, if $f(k)=(ak+b)\bi{2k}k^2\bi{3k}k/m^k$ with $m>0$, then 
\begin{align*}\sum_{k=0}^\infty f'(k)=&\ \sum_{k=0}^\infty \l((ak+b)(2H_{2k}+3H_{3k}-5H_k)+a\r)\f{\bi{2k}k^2\bi{3k}k}{m^k}
\\&\ -(\log m)\sum_{k=0}^\infty(ak+b)\f{\bi{2k}k^2\bi{3k}k}{m^k}
\end{align*}
vanishes. For a Ramanujan-type series $\sum_{k=0}^\infty f(k)$ of the form \eqref{R},
there are some observations about $\sum_{k=0}^\infty f''(k)$ and $\sum_{k=0}^\infty f'''(k)$
(cf. the author's posting \cite{MO} and the preprint \cite{CG}). They will be proved  
in a forthcoming paper by H. Cohen, the author and Y. Zhou.

Motivated by the known identity
 $$\sum_{k=1}^\infty\f1{k^2\binom{2k}k}=\f{\pi^2}{18}=\f{\zeta(2)}3,$$
 D. Chen \cite{Chen} guessed the value of $\sum_{k=1}^\infty f^{(m)}(k)$ for $m\ls5$ in terms of
 Riemann's zeta function, where
 $f(k)=\Gamma(k+1)^2/(k^2\Gamma(2k+1))$.

 In Sections 4 and 5, we propose many new conjectural series with summands involving higher-order derivatives.

 As all the conjectural series in this paper converge fast, the conjectural series identities
 can be easily checked numerically via {\tt Mathematica}.

 Now we introduce some notations throughout this paper.

 For and odd prime $p$ and integer $a$, as usual $(\f ap)$ denotes the Legendre symbol.
 For any non-square integer $d\eq0,1\pmod4$ and real number $s\gs1$, we define
 $$L_d(s):=L\left(s,\l(\frac{d}{\cdot}\r)\right)=\sum_{n=1}^\infty\f{(\f dn)}{n^s}$$
 with $(\f dk)$ the Kronecker symbol. In particular, for $s\gs1$ we have
$$L_{-3}(s)=\sum_{n=1}^\infty\frac{(\f n3)}{n^s}
=\sum_{k=0}^\infty\left(\frac1{(3k+1)^s}-\frac1{(3k+2)^s}\right)$$
and
$$L_{-8}(2)=\sum_{n=1}^\infty\f{(\f{-8}n)}{n^s}=\sum_{k=0}^\infty\f{(-1)^{k(k-1)/2}}{(2k+1)^s}.$$
The Dirichlet beta function is given by
$$\beta(s)=\sum_{k=0}^\infty\f{(-1)^k}{(2k+1)^s}\quad\t{with}\ s\gs1.$$
Note that $G=\beta(2)$ is called {\it Catalan's constant}.
Some of our conjectural identities involve the multiple zeta value
$$\zeta(5,3):=\sum_{k_1>k_2>0}\frac1{k_1^5k_2^3}=0.0377\ldots.$$

\section{A new type of series for $1/\pi$ involving generalized central trinomial coefficients}
\setcounter{equation}{0}
 \setcounter{conjecture}{0}
 \setcounter{theorem}{0}
 \setcounter{proposition}{0}

 Now we introduce series for $1/\pi$ involving generalized central trinomial coefficients
 of the following new type:

\ \  {\tt Type X}. $\sum_{k=0}^\infty \f{a+dk}{m^k}\bi{2k}kT_{k}(b,c)T_{2k}(b_*,c_*)$,

\noindent where $a\in\Z$ and $b,c,b_*,c_*,d,m\in\Z\sm\{0\}$ with $b^2\not=4c$ and $b_*^2\not=4c_*$.

 \begin{conjecture} We have the series identities
 \[\label{X1}\sum_{k=0}^\infty\f{24k+5}{76^{2k}}\bi{2k}kT_k(1,-12)T_{2k}(8,-3)
 =\f{\sqrt{38(120+73\sqrt3)}}{6\pi},
 \tag{X1}\]
  \[\label{X2}\sum_{k=0}^\infty\f{585k+136}{(-85^2)^k}\bi{2k}kT_k(8,-9)T_{2k}(7,-9)
 =\f{85\sqrt{255}}{6\pi},
 \tag{X2}\]
 \[\label{X3}\sum_{k=0}^\infty\f{16k+3}{(-202^2)^k}
 \bi{2k}kT_k(19,-20)T_{2k}(9,-5)=\f{43\sqrt{101}}{75\pi}.\tag{X3}
 \]
 \end{conjecture}
 \begin{remark} \eqref{X1} and \eqref{X2} were found on March 18, 2026, and
 \eqref{X3} was discovered on March 20, 2026.
 \end{remark}

 Below are related conjectures on congruences.

\begin{conjecture} [2026-03-18] \label{Conj2.2}
 {\rm (i)} For any integer $n>1$, we have
 \begin{equation}\f1{n\bi{2n}n}\sum_{k=0}^{n-1}(24k+5)76^{2(n-1-k)}\bi{2k}kT_k(1,-12)T_{2k}(8,-3)
 \in\Z^+,
 \end{equation}
 and this number is odd if and only if $n\in\{2^a+1:\ a\in\N\}$.

{\rm (ii)} For any prime $p>3$ with $p\not=19$, we have
\begin{equation}
 \begin{aligned}
 &\ \sum_{k=0}^{p-1}\f{\bi{2k}k T_k(1,-12) T_{2k} (8,-3)}{76^{2k}}
 \\\eq&\ \begin{cases}(\f p{19})(4x^2-2p)\pmod{p^2}&\t{if}\ p=x^2+4y^2\ (x,y\in\Z),
 \\0\pmod{p^2}.\end{cases}
 \end{aligned}
 \end{equation}
 \end{conjecture}
 \begin{remark}  Actually, \eqref{X1} was motivated by this conjecture.
 \end{remark}

 \begin{conjecture} [2026-03-18]\label{Conj2.3}
 {\rm (i)} For any $n\in\Z^+$, we have
 $$\f1{2n\bi{2n}n}\sum_{k=0}^{n-1}(-1)^k(585k+136)85^{2(n-1-k)}\bi{2k}kT_k(8,-9)T_{2k}(7,-9)\in\Z^+,$$
 and this number is odd if and only if $n\in\{2^a:\ a\in\Z^+\}$.

 {\rm (ii)} For any odd prime $p\not=5,17$, we have
 \begin{equation}\begin{aligned}&\ \sum_{k=0}^{p-1}\f{585k+136}{(-7225)^k}
 \bi{2k}kT_k(8,-9)T_{2k}(7,-9)
 \\\eq&\ p\l(\f p{17}\r)\l(75\l(\f p{15}\r)+\f{1053-16(\f p5)}{17}\r)\pmod{p^2}.
 \end{aligned}
 \end{equation}

 {\rm (iii)} For any prime $p>5$ with $p\not=17$, we have
 \begin{equation}\begin{aligned}&\ \l(\f p{17}\r)\sum_{k=0}^{p-1}\f{\bi{2k}kT_k(8,-9)T_{2k}(7,-9)}{(-7225)^k}
 \\\eq&\ \begin{cases}4x^2-2p\pmod{p^2}&\t{if}\ p\eq1,4\pmod{15}\ \&\ p=x^2+15y^2,
 \\2p-12x^2\pmod{p^2}&\t{if}\ p\eq 2,8\pmod{15}\ \&\ p=3x^2+5y^2,
 \\0\pmod{p^2}&\t{if}\ (\f {-15}p)=-1,
 \end{cases}\end{aligned}
 \end{equation}
 where $x$ and $y$ are integers.
 \end{conjecture}
 \begin{remark} \eqref{X2} was motivated by this conjecture.
 \end{remark}

 \begin{conjecture} [2026-03-20] \label{Conj2.4}
 {\rm (i)} For any $n\in\Z^+$, the number
 $$\f3{n\bi{2n-1}{n-1}}\sum_{k=0}^{n-1}(-1)^k(16k+3)202^{2(n-1-k)}\bi{2k}kT_k(19,-20)T_{2k}(9,-5)$$
 is always a positive odd integer.

 {\rm (ii)} For any odd prime $p\not=5,101$, we have
 \begin{equation}\begin{aligned}&\ \sum_{k=0}^{p-1}\f{16k+3}{(-202^2)^k}
 \bi{2k}kT_k(19,-20)T_{2k}(9,-5)
 \\\eq&\ \f p{7575}\l(\f{p}{101}\r)\l(16856\l(\f{-1}p\r)+6000\l(\f p5\r)-131\r)
 \pmod{p^2}.\end{aligned}\end{equation}

 {\rm (iii)} For any odd prime $p\not=5,101$, we have
 \begin{equation}\begin{aligned} &\l(\f p{101}\r)\sum_{k=0}^{p-1}
 \f{\bi{2k}kT_k(19,-20)T_{2k}(9,-5)}{(-202^2)^k}
 \\\eq& \begin{cases}4x^2-2p\pmod{p^2}&\t{if}\ p\eq1,9\pmod{20}\ \&\ p=x^2+y^2\ (5\nmid x),
 \\4xy\pmod{p^2}&\t{if}\ p\eq13,17\pmod{20}\ \&\ p=x^2+y^2\ (5\mid x-y),
 \\0\pmod{p^2},\end{cases}
 \end{aligned}
 \end{equation}
 where $x$ and $y$ are integers.
 \end{conjecture}
 \begin{remark}  \eqref{X3} was motivated by this conjecture.
 \end{remark}

Our following six conjectures are similar to Conjectures \ref{Conj2.2}-\ref{Conj2.4}.

\begin{conjecture} [2026-03-14] For any odd prime $p\not=5$, we have
 \begin{equation}
 \begin{aligned}
 &\ \sum_{k=0}^{p-1} \f{\bi{2k}kT_k(3,-4)T_{2k}(1,-1)}{(-100)^k}
 \\\eq&\ \begin{cases}4x^2-2p\pmod{p^2}&\t{if}\ p\eq 1,9\pmod{20}\ \&\ p=x^2+5y^2,
 \\2p-2x^2\pmod{p^2}&\t{if}\ p\eq3,7\pmod{20}\ \&\ 2p=x^2+5y^2,
 \\0\pmod{p^2}&\t{if}\ (\f{-5}p)=-1,
 \end{cases}
 \end{aligned}
 \end{equation}
 where $x$ and $y$ are integers.
 \end{conjecture}

 \begin{conjecture} [2026-03-14] For any prime $p>5$, we have
 \begin{equation}
 \begin{aligned}
 &\ \l(\f p5\r)\sum_{k=0}^{p-1}\f{\bi{2k}kT_k(5,4)T_{2k}(3,1)}{100^k}
 \\\eq&\ \begin{cases}4x^2-2p\pmod{p^2}&\t{if}\ p\eq 1,19\pmod{24}\ \&\ p=x^2+6y^2,
 \\8x^2-2p\pmod{p^2}&\t{if}\ p\eq5,23\pmod{24}\ \&\ p=2x^2+3y^2,
 \\0\pmod{p^2}&\t{if}\ (\f{-6}p)=-1,
 \end{cases}
 \end{aligned}
 \end{equation}
 where $x$ and $y$ are integers.
 \end{conjecture}

 \begin{conjecture} [2026-03-15] Let $p>3$ be a prime. Then
 \begin{equation}\begin{aligned}&\ \l(\f{-1}p\r)\sum_{k=0}^{p-1}\f{\bi{2k}kT_k(10,9)T_{2k}(2,9)}{1024^k}
 \\\eq&\begin{cases}4x^2-2p\pmod{p^2}&\t{if}\ p\eq1,7\pmod{24}\ \&\ p=x^2+6y^2\,(x,y\in\Z),
 \\8x^2-2p\pmod{p^2}&\t{if}\ p\eq5,11\pmod{24}\ \&\ p=2x^2+3y^2\,(x,y\in\Z),
 \\0\pmod{p^2}&\t{if}\ (\f{-6}p)=-1.
 \end{cases}\end{aligned}
 \end{equation}
 If $(\f{-6}p)=1$, then
 \begin{equation}\sum_{k=0}^{p-1}\f{144k+43}{1024^k}\bi{2k}kT_k(10,9)T_{2k}(2,9)\eq 6p\l(\f{-1}p\r)\pmod{p^2}.
 \end{equation}
 \end{conjecture}

  \begin{conjecture} [2026-03-19] For any odd prime $p\not=5,13$, we have
 \begin{equation}
 \begin{aligned}&\l(\f p{13}\r)\sum_{k=0}^{p-1}\f{\bi{2k}kT_k(17,16)T_{2k}(9,4)}{130^{2k}}
 \\\eq&\begin{cases}4x^2-2p\pmod{p^2}&\t{if}\ (\f{-2}p)=(\f 5p)=1\ \&\ p=x^2+10y^2,
 \\2p-8x^2\pmod{p^2}&\t{if}\ (\f{-2}p)=(\f 5p)=-1\ \&\ p=2x^2+5y^2,
 \\0\pmod{p^2}&\t{if}\ (\f{-10}p)=-1,
 \end{cases}
 \end{aligned}
 \end{equation}
 where $x$ and $y$ are integers.
 \end{conjecture}

 \begin{conjecture} [2026-03-17] For any odd prime $p\not=37$, we have
 \begin{equation}
 \begin{aligned}&\l(\f p{37}\r)\sum_{k=0}^{p-1}\f{\bi{2k}kT_k(13,12)T_{2k}(7,3)}{74^{2k}}
 \\\eq&\begin{cases}4x^2-2p\pmod{p^2}&\t{if}\ p\eq1\pmod{12}\ \&\ p=x^2+y^2\ (3\nmid x),
 \\4xy\pmod{p^2}&\t{if}\ p\eq5\pmod{12}\ \&\ p=x^2+y^2\ (3\mid x-y),
 \\0\pmod{p^2}&\t{if}\ p\eq3\pmod4,
 \end{cases}
 \end{aligned}
 \end{equation}
 where $x$ and $y$ are integers.
 Also, for any prime $p\eq3\pmod4$ we have
 $$\sum_{k=0}^{p-1}\f{k\bi{2k}kT_k(13,12)T_{2k}(7,3)}{74^{2k}}\eq0\pmod p.$$
 \end{conjecture}

\begin{conjecture} [2026-03-27] For any odd prime $p\not=13$, we have
 \begin{equation}
 \begin{aligned}&\l(\f p{13}\r)\sum_{k=0}^{p-1}\f{\bi{2k}kT_k(25,24)T_{2k}(14,-3)}{104^{2k}}
 \\\eq&\begin{cases}4x^2-2p\pmod{p^2}&\t{if}\ p\eq1\pmod{12}\ \&\ p=x^2+y^2\ (3\nmid x),
 \\4xy\pmod{p^2}&\t{if}\ p\eq5\pmod{12}\ \&\ p=x^2+y^2\ (3\mid x-y),
 \\0\pmod{p^2}&\t{if}\ p\eq3\pmod4,
 \end{cases}
 \end{aligned}
 \end{equation}
 where $x$ and $y$ are integers.
 Also, for any prime $p>3$ with $p\eq3\pmod4$ with $p\not=19$, we have
 \begin{equation}\sum_{k=0}^{p-1}\f{k\bi{2k}kT_k(25,24)T_{2k}(14,-3)}{104^{2k}}\eq0\pmod p.
 \end{equation}
 \end{conjecture}

 \begin{conjecture} [2026-04-01]  {\rm (i)} For any integer $n>1$, the number
 $$\f1{3n\bi{2n-1}{n-1}}\sum_{k=0}^{n-1}(76k+35)2700^{n-1-k}\bi{2k}kT_k(31,-32)T_{2k}(7,1)$$
 is always an odd integer.

 {\rm (ii)} For any prime $p>5$, we have
 \begin{equation}
 \begin{aligned}&\ \l(\f p5\r)\sum_{k=0}^{p-1}\f{\bi{2k}kT_k(31,-32)T_{2k}(7,1)}{2700^k}
 \\\eq&\ \begin{cases}4x^2-2p\pmod{p^2}&\t{if}\ p\eq1\pmod{12}\ \&\ p=x^2+y^2\ (3\nmid x),
 \\4xy\pmod{p^2}&\t{if}\ p\eq5\pmod{12}\ \&\ p=x^2+y^2\ (3\mid x-y),
 \\0\pmod{p^2}&\t{if}\ p\eq3\pmod4,\end{cases}
 \end{aligned}\end{equation}
 where $x$ and $y$ are integers.
 \end{conjecture}

 \begin{conjecture} [2026-03-29]  {\rm (i)} For any $n\in\Z^+$, the number
 $$\f1{n\bi{2n-1}{n-1}}\sum_{k=0}^{n-1}(96k+25)22^{2(n-1-k)}\bi{2k}kT_k(37,36)T_{2k}(5,9)$$
 is always an odd integer.

 {\rm (ii)} Let $p$ be an odd prime with $p\not=11$. Then
 \begin{equation}
 \begin{aligned}&\ \sum_{k=0}^{p-1}\f{\bi{2k}kT_k(37,36)T_{2k}(5,9)}{22^{2k}}
 \\\eq&\ \begin{cases}(\f{33}p)(4x^2-2p)\pmod{p^2}&\t{if}\ (\f{-2}p)=1\ \&\ p=x^2+2y^2\ (x,y\in\Z),
 \\0\pmod{p^2}&\t{if}\ (\f{-2}p)=-1,\ \t{i.e.},\ p\eq5,7\pmod8.\end{cases}
 \end{aligned}\end{equation}
 \end{conjecture}

 \begin{conjecture} [2026-03-29] For any prime $p>3$ with $p\eq3\pmod4$, we have
 \begin{equation}\sum_{k=0}^{p-1}\f{k\bi{2k}kT_k(3,2)T_{2k}(4,-2)}{(-288)^k}
 \eq \sum_{k=0}^{p-1}\f{\bi{2k}kT_k(3,2)T_{2k}(4,-2)}{(-288)^k}\eq0\pmod p.
 \end{equation}
 \end{conjecture}

 \section{Conjectural series involving harmonic numbers}
 \setcounter{equation}{0}
 \setcounter{conjecture}{0}
 \setcounter{theorem}{0}
 \setcounter{proposition}{0}

  Below we pose some new conjectures on series involving (generalized) harmonic numbers.

 \begin{conjecture} [2026-01-25] We have
 \begin{equation}\sum_{k=0}^\infty\f{\bi{2k}k^2\bi{3k}k}{(-192)^{k}}
 \l((5k+1)(H_{2k}-H_k)+\f{3k+1}{2k+1}\r)=\f{4\log3}{\sqrt3\,\pi}
 \end{equation}
 and
 \begin{equation}\sum_{k=0}^\infty\f{\bi{2k}k^2\bi{3k}k}{(-192)^{k}}
 \l((5k+1)(H_{3k}-H_k)+\f{4k+3}{6k+3}\r)=\f4{3\sqrt3\,\pi}\log\f{64}3.
 \end{equation}
 \end{conjecture}
 \begin{remark} This is motivated by the Ramanujan series
 $$\sum_{k=0}^\infty(5k+1)\f{\bi{2k}k^2\bi{3k}k}{(-192)^{k}}=\f{4}{\sqrt3\,\pi}.$$
 \end{remark}

 \begin{conjecture} [2026-01-25] We have
 \begin{equation}\sum_{k=0}^\infty\f{\bi{2k}k^2\bi{3k}k}{(-12)^{3k}}
 \l((51k+7)(H_{2k}+H_k)-\f{3(13k+7)}{2k+1}\r)=-36\sqrt3\f{\log3}{\pi}
 \end{equation}
 and
 \begin{equation}\sum_{k=0}^\infty\f{\bi{2k}k^2\bi{3k}k}{(-12)^{3k}}
 \l((51k+7)(3H_{3k}-7H_k)-\f{180k+93}{2k+1}\r)=36\sqrt3\f{\log108}{\pi}
 \end{equation}
 \end{conjecture}
 \begin{remark} This is motivated by the known series
 $$\sum_{k=0}^\infty(51k+7)\f{\bi{2k}k^2\bi{3k}k}{(-12)^{3k}}=\f{12\sqrt3}{\pi}.$$
 \end{remark}

 \begin{conjecture} [2026-01-25] We have
 \begin{equation}\sum_{k=0}^\infty\f{\bi{2k}k^2\bi{3k}k}{(-8640)^{k}}
 \l((9k+1)(H_{2k}-2H_k)-\f{27k+13}{6k+3}\r)=\f{48\log 3}{\sqrt{15}\,\pi}
 \end{equation}
 and
 \begin{equation}\sum_{k=0}^\infty\f{\bi{2k}k^2\bi{3k}k}{(-8640)^{k}}
 \l((9k+1)(3H_{3k}-H_k)+\f{1}{6k+3}\r)=\f{12}{\sqrt{15}\,\pi}\log\f{320}{243}.
 \end{equation}
 \end{conjecture}
 \begin{remark} This is motivated by the known series
 $$\sum_{k=0}^\infty(9k+1)\f{\bi{2k}k^2\bi{3k}k}{(-8640)^{k}}=\f{4\sqrt{15}}{5\pi}.$$
 \end{remark}

 \begin{conjecture} [2026-01-26]
We have
\begin{equation}
\sum_{k=0}^\infty\f{\bi{2k}k^2\bi{4k}{2k}}
{(-12288)^{k}}\l((28k+3)(H_{2k}-H_k)+\f{20k+9}{4k+2}\r)=\f{4\log432}{\sqrt3\,\pi},
\end{equation}
\begin{equation}
\sum_{k=0}^\infty\f{\bi{2k}k^2\bi{4k}{2k}}
{(-12288)^{k}}\l((28k+3)(H_{4k}-H_{2k})+\f{8k+5}{4k+2}\r)=\f{8}{\sqrt3\,\pi}\log\f{16}3
\end{equation}
\end{conjecture}
\begin{remark} This is motivated by the known series
$$\sum_{k=0}^\infty(28k+3)\f{\bi{2k}k^2\bi{4k}{2k}}{(-12288)^{k}}=\f{16}{\sqrt3\,\pi}.$$
\end{remark}

\begin{conjecture} [2026-01-26] We have
\begin{equation}
\sum_{k=0}^\infty\f{\bi{2k}k^2\bi{4k}{2k}}
{(-12288)^{k}}\l((28k+3)(36H_{4k}^{(2)}-9H_{2k}^{(2)}-H_k^{(2)})+\f{12}{4k+1}\r)=\f{16\pi}{3\sqrt3}.
\end{equation}
\end{conjecture}
\begin{remark} This is quite similar to \cite[Conjecture 41]{harmonic}.
\end{remark}

\begin{conjecture} [2026-01-26]
We have
\begin{equation}
\sum_{k=0}^\infty\f{\bi{2k}k\bi{3k}k\bi{6k}{3k}}
{(-15)^{3k}}\l((63k+8)(3H_{2k}-5H_k)+\f{198k+83}{4k+2}\r)=30\sqrt{15}\f{\log3}{\pi}.
\end{equation}
\end{conjecture}
\begin{remark} This is motivated by the known series
$$\sum_{k=0}^\infty(63k+8)\f{\bi{2k}k\bi{3k}k\bi{6k}{3k}}{(-15)^{3k}}=\f{5\sqrt{15}}{\pi}.$$
\end{remark}

\begin{conjecture} [2026-01-25]
We have
\begin{equation}
\sum_{k=0}^\infty\f{\bi{2k}k\bi{3k}k\bi{6k}{3k}}
{(-96)^{3k}}\l((342k+25)(3H_{2k}-7H_k)+\f{1122k+571}{2k+1}\r)=64\sqrt{6}\f{\log(2^73^6)}{\pi}.
\end{equation}
\end{conjecture}
\begin{remark} This is motivated by the known series
$$\sum_{k=0}^\infty(342k+25)\f{\bi{2k}k\bi{3k}k\bi{6k}{3k}}{(-96)^{3k}}=\f{32\sqrt{6}}{\pi}.$$
\end{remark}

 \begin{conjecture} [2026-01-24] We have
 \begin{equation}\sum_{k=1}^\infty\f{(21k-8)(2H_{3k-1}-H_{k-1})-\f{4(27k-11)}{3(3k-1)}}{k^3\bi{2k}k^3}
 =3\zeta(3).
 \end{equation}
 \end{conjecture}
 \begin{remark} This is motivated by the Zeilberger series (cf. \cite{Zeil})
 $$\sum_{k=1}^\infty\f{21k-8}{k^3\bi{2k}k^3}=\f{\pi^2}6.$$
 \end{remark}

 \begin{conjecture} [2026-01-26] We have
 \begin{equation}\sum_{k=1}^\infty\f{(-27)^{k-1}((15k-4)(9H_{3k-1}^{(2)}-5H_{k-1}^{(2)})-18/(3k-1))}
 {k^3\bi{2k}k^2\bi{3k}k}=9L_{-3}(4).
 \end{equation}
 \end{conjecture}
 \begin{remark} This is motivated by the identity
 $$\sum_{k=1}^\infty\f{(15k-4)(-27)^{k-1}}{k^3\bi{2k}k^2\bi{3k}k}=L_{-3}(2)$$
 conjectured by the author \cite{S11} and later confirmed by J. Guillera and M. Rogers \cite{GR}.
 \end{remark}

 \begin{conjecture} [2026-03-09] We have
 \begin{equation}\sum_{k=1}^\infty\f{H_{k-1}^{(2)}}{k^2\bi{2k}k}\l(H_{2k}+26H_k-\f 8k\r)
 -\f{34\pi^2\zeta(3)-287\zeta(5)}{54}.
 \end{equation}
 \end{conjecture}
 \begin{remark} This is motivated by the known identity
 $$\sum_{k=1}^\infty\f{H_{k-1}^{(2)}}{k^2\bi{2k}k}=\f{\pi^4}{1944}$$
 which follows from the fact that
 $$\arcsin^4\f x2=\f32\sum_{k=1}^\infty\f{H_{k-1}^{(2)}}{k^2\bi{2k}k}x^{2k}\ \ \t{for all}\ -2\ls x\ls 2$$
 (cf. \cite{BC}).
 \end{remark}

\begin{conjecture} [2026-04-12] We have
\begin{equation}\label{21}\begin{aligned}&\ \sum_{k=1}^\infty\f{256^k\l(14(21k^3-22k^2+8k-1)(H_{2k-1}-H_{k-1})-63k^2+44k-8\r)}
{k^7\bi{2k}k^7}\\&\qquad =\f{93}2\zeta(5)+\pi^4\log2.
\end{aligned}
\end{equation}
\end{conjecture}
\begin{remark} This is motivated by the identity
$$\sum_{k=1}^\infty \f{(21k^3-22k^2+8k-1)256^k}{k^7\bi{2k}k^7}=\f{\pi^4}8$$
conjectured by Guillera in 2003 and confirmed by Au \cite[Example XI]{Au} in 2025.
Actually, \eqref{21} is equivalent to the identity
$$\sum_{k=1}^\infty f'(k)=-\f{93}2\zeta(5),$$
where
$$f(x)=2(21x^3-22x^2+8x-1)256^{x-1}\f{\Gamma(x)^{14}}{\Gamma(2x)}\quad \t{for}\ x>0.$$
\end{remark}

\section{Conjectural series involving higher-order derivatives (I)}
\setcounter{equation}{0}
 \setcounter{conjecture}{0}
 \setcounter{theorem}{0}
 \setcounter{proposition}{0}

 In this section, we pose many conjectures on series with summands
 containing no nontrivial powers and involving higher-order derivatives.

\begin{conjecture} [2026-03-05] For $x>0$, set
$$f_1(x)=\frac{\Gamma(x)^2}{2\Gamma(2x)}.$$
Then
\begin{align}\sum_{k=1}^\infty f_1'(k)&=-\frac32L_{-3}(2),
\\\sum_{k=1}^\infty f_1''(k)&=3L_{-3}(3)=\frac{4\pi^3}{27\sqrt3},
\\\sum_{k=1}^\infty f_1'''(k)&=\frac34\left(2\pi^2L_{-3}(2)-27L_{-3}(4)\right).
\end{align}
\end{conjecture}
\begin{remark} It is well known that
$$\sum_{k=1}^\infty f_1(k)=\sum_{k=1}^\infty\frac1{k\binom{2k}k}=\frac{\pi}{3\sqrt3}=L_{-3}(1).$$
\end{remark}

\begin{conjecture} [2026-03-09] For $x>0$, define
$$f_2(x)=\f{(\zeta(2)-\sum_{n=0}^\infty (x+n)^{-2})\Gamma(x)^2}{2x\Gamma(2x)}.$$
Then
\begin{align}\sum_{k=1}^\infty f_2'(k)&=\f{29}{27}\zeta(5),
\\\sum_{k=1}^\infty f_2''(k)&=-\f{4115}{648}\zeta(6)-\f 83\zeta(3)^2,
\\\sum_{k=1}^\infty f_2'''(k)&=24\zeta(3)\zeta(4)+\f1{18}(1069\zeta(7)-116\zeta(2)\zeta(5)),
\\\sum_{k=1}^\infty f_2^{(4)}(k)&=-\f{359939}{648}\zeta(8)-128\zeta(3)\zeta(5)
+32\zeta(2)\zeta(3)^2+\f{224}9\zeta(5,3),
\end{align}
and
\begin{equation}\begin{aligned}\sum_{k=1}^\infty f_2^{(5)}(k)&=\f{459385}{81}\zeta(9)-\f{10690}9\zeta(2)\zeta(7)
-\f{640}3\zeta(3)^3
\\&\ \ \ +1990\zeta(4)\zeta(5)-40\zeta(3)\zeta(6).
\end{aligned}
\end{equation}
\end{conjecture}
\begin{remark} Note that
$$\sum_{k=1}^\infty f_2(k)=\sum_{k=1}^\infty\f{H_{k-1}^{(2)}}{k^2\bi{2k}k}=\f{\pi^4}{1944}$$
as mentioned in Remark 3.10.
\end{remark}

\begin{conjecture} [2026-03-05] For $x>0$, set
$$ f_3(x)=\frac{\Gamma(x)^2}{2x^3\Gamma(2x)}.$$
Then
\begin{align}\sum_{k=1}^\infty f_3''(k)&=\frac{3073}{216}\zeta(6),
\\\sum_{k=1}^\infty f_3^{(3)}(k)&=\frac{176\zeta(2)\zeta(5)-1423\zeta(7)}{12},
\\\sum_{k=1}^\infty f_3^{(4)}(k)&=\frac{752537\zeta(8)+25344\zeta(5,3)}{1080},
\\\sum_{k=1}^\infty f_4^{(5)}(k)&=660\zeta(4)\zeta(5)+\frac{7115}3\zeta(2)\zeta(7)-\frac{283010}{27}\zeta(9).
\end{align}
\end{conjecture}
\begin{remark} It is well known that
$$\sum_{k=1}^\infty f_3(k)=\sum_{k=1}^\infty\frac1{k^4\binom{2k}k}=\frac{17}{36}\zeta(4).$$
Note also that $\sum_{k=1}^\infty f_3'(k)=-22\zeta(5)/9$ which is equivalent to the identity
$$\sum_{k=1}^\infty\f{H_{2k}-H_k+2/k}{k^4\bi{2k}k}=\f{11}9\zeta(5)$$
conjectured by the author \cite{S15} and later confirmed by Ablinger \cite{Abl}.
\end{remark}

\begin{conjecture} [2026-03-05] For $x>1/2$, define
$$f_4(x):=\frac{e^{\pi ix}(28x^2-18x+3)\Gamma(x)^3}{3x^2(2x-1)^3\Gamma(3x)}.$$
Then
\begin{align}\sum_{k=1}^\infty f_4'(k)&=20\zeta(5)-i\frac{\pi^5}{45},
\\\sum_{k=1}^\infty f_4''(k)&=198\zeta(6)-40i\pi\zeta(5),
\\\sum_{k=1}^\infty f_4'''(k)&=36(91\zeta(7)-20\zeta(2)\zeta(5))-\f{212}{315}\pi^7i,
\end{align}
\begin{equation}\sum_{k=1}^\infty f_4^{(4)}(k)=16\pi i(819\zeta(7)-20\pi^2\zeta(5))-96(243\zeta(8)+8\zeta(5,3)),
\end{equation}
and
\begin{equation}\begin{aligned}\sum_{k=1}^\infty f_4^{(5)}(k)=&\ 4320(208\zeta(9)-91\zeta(2)\zeta(7)+10\zeta(4)\zeta(5))
\\&\ -\pi i(665\zeta(8)+16\zeta(5,3)).
\end{aligned}
\end{equation}
\end{conjecture}
\begin{remark} Note that we have the identity
$$ \sum_{k=1}^\infty f_4(k)=\sum_{k=1}^\infty\frac{(-1)^k(28k^2-18k+3)}{(2k-1)^3k^4\binom{2k}k\binom{3k}k}=-\frac{\pi^4}{45},$$
as conjectured by the author \cite[(4.12)]{S24} and confirmed by K. C. Au \cite{Au26} recently.
\end{remark}

\begin{conjecture} [2026-03-07] For $x>0$, define
$$f_5(x)=\f{e^{\pi i x}S(x)\Gamma(x)^3}{3x^3(2x-1)^4\Gamma(3x)},$$
where
 $$S(x)=560x^4-640x^3+480x^2-136x+17.$$
 Then
 \begin{align}\label{5'}\sum_{k=1}^\infty f_5'(k)&=392\zeta(3)^2-190\zeta(6)+4\pi i\l(45\zeta(5)-\f{14}3\pi^2\zeta(3)\r),
 \\\sum_{k=1}^\infty f_5''(k)&=240(16\zeta(7)-9\zeta(2)\zeta(5)-14\zeta(3)\zeta(4))+4\pi i(196\zeta(3)^2-95\zeta(6)),
 \end{align}
 \begin{equation}
 \begin{aligned}\sum_{k=1}^\infty f_5'''(k)=&\ \f{18}5\l(1711\zeta(8)+17360\zeta(3)\zeta(5)
 +2992\zeta(5,3)-3920\zeta(2)\zeta(3)^2\r)
 \\&\ +480\pi i\l(24\zeta(7)-9\zeta(2)\zeta(5)-28\zeta(3)\zeta(4)\r),
 \end{aligned}\end{equation}
 \begin{equation}
 \begin{aligned}&\ \f{\sum_{k=1}^\infty f_5^{(4)}(k)}{64}
 \\=&\ 392\zeta(3)^3-3695\zeta(9)-4320\zeta(2)\zeta(7)-1575\zeta(3)\zeta(6)+1215\zeta(4)\zeta(5)
 \\&\ +\f{\pi i}{40}\l(5899\zeta(8)-23520\zeta(2)\zeta(3)^2+156240\zeta(3)\zeta(5)+26928\zeta(5,3)\r),
 \end{aligned}\end{equation}
 and
 \begin{equation}\f{\sum_{k=1}^\infty \mathrm{Im}(f_5^{(5)}(k))}{320\pi}=392\zeta(3)^3-3695\zeta(9)
 -2880\zeta(2)\zeta(7)-4662\zeta(3)\zeta(6).
 \end{equation}
 \end{conjecture}
 \begin{remark} This is motivated by the identity
 $$\sum_{k=1}^\infty f_5(k)=\sum_{k=1}^\infty \f{(-1)^kS(k)}{k^5(2k-1)^4\bi{2k}k\bi{3k}k}=180\zeta(5)-\f{56}3\pi^2\zeta(3)$$
 conjectured by the author in Jan. 2025 (cf. \cite{zeta(5),II}) and confirmed by K. C. Au via the WZ-seed method (cf. Au's answer in \cite{II}). We observe that \eqref{5'} has the following equivalent form:
 \begin{equation*}\begin{aligned} &\sum_{k=1}^\infty\f{(-1)^k(S(k)(H_{3k-1}-H_{k-1})+40T(k)+\f{50k-17}{k(2k-1)})}{k^5(2k-1)^4\bi{2k}k\bi{3k}k}
=\f23\l(95\zeta(6)-196\zeta(3)^2\r),
 \end{aligned}\end{equation*}
 where $T(k)=14k^3-12k^2+11k-3.$
 \end{remark}

\begin{conjecture} For $x>1/4$, define
$$f_{6}(x)=\f{P(x)\Gamma(4x-1)\Gamma(x)^{16}}{256\Gamma(2x)^{10}},$$
where
$$P(x)=5460x^4-8341x^3+4864x^2-1280x+128.$$
Then
\begin{align}\label{f6'}\sum_{k=1}^\infty f_{6}'(k)&=-48\zeta(5),
\\\sum_{k=1}^\infty f_{6}''(k)&=356\zeta(6)=\f{356}{945}\pi^6,
\\\sum_{k=1}^\infty f_{6}'''(k)&=288(4\zeta(2)\zeta(5)-15\zeta(7)),
\\\sum_{k=1}^\infty f_{6}^{(4)}(k)&=12512\zeta(8)=\f{6256}{4725}\pi^8,
\\\sum_{k=1}^\infty f_{6}^{(5)}(k)&=5760(60\zeta(2)\zeta(7)-48\zeta(4)\zeta(5)-49\zeta(9)),
\\\sum_{k=1}^\infty f_{6}^{(6)}(k)&=1728\l(1440\zeta(5)^2-1999\zeta(10)\r),
\end{align}
and
\begin{equation}\f{\sum_{k=1}^\infty f_{6}^{(7)}(k)}{241920}
=885\zeta(11)+196\zeta(2)\zeta(9)-720\zeta(4)\zeta(7)-336\zeta(5)\zeta(6).
\end{equation}
\end{conjecture}
\begin{remark} Note that
$$\sum_{k=1}^\infty f_{6}(k)=\sum_{k=1}^\infty \f{P(k)\bi{4k}{2k}}{(4k-1)k^7\bi{2k}k^8}$$
was conjectured to have the value $\pi^4/15=6\zeta(4)$ by D. Chen \cite{Chen25}.
We observe that \eqref{f6'} is equivalent to the identity
$$\sum_{k=1}^\infty\f{(P(k)(H_{4k-1}-5H_{2k-1}+4H_{k-1})+(4k-1)\f{Q(k)}4-\f3{4k-1})\bi{4k}{2k}}
{(4k-1)k^7\bi{2k}k^8}=-12\zeta(5),$$
where $Q(k)=4095k^2-3488k+780.$

We have many other similar conjectures which will not be listed here.
\end{remark}

\section{Conjectural series involving higher-order derivatives (II)}
\setcounter{equation}{0}
 \setcounter{conjecture}{0}
 \setcounter{theorem}{0}
 \setcounter{proposition}{0}

 In this section, we pose many conjectures on series with summands
 containing nontrivial powers and involving higher-order derivatives.

 \begin{conjecture} [2026-04-10] For $x>-1/6$, define
 $$g_1(x)=\f{(10x-1)\Gamma(6x+1)}{(2x+1)512^x\Gamma(x+1)\Gamma(2x+1)\Gamma(3x+1)}.$$
 Then
 \begin{equation}
 \sum_{k=0}^\infty g_1'(k)=4\sqrt2\,\pi
 \ \ \t{and}\ \ \sum_{k=0}^\infty g_1''(k)=-256L_{-8}(2).
 \end{equation}
 \end{conjecture}
 \begin{remark} Note that
 $$\sum_{k=0}^\infty g_1(k)=\sum_{k=0}^\infty \f{(10k-1)\bi{6k}{3k}\bi{3k}k}{(2k+1)512^k}=0$$
 as conjectured by the author \cite[Remark 4.6]{S23} and later confirmed by Au \cite{Au}.
 \end{remark}

 \begin{conjecture} [2026-04-10] For $x>-1/6$, define
 $$g_2(x)=\f{(74x+7)\Gamma(6x+1)}{(2x+1)4096^x\Gamma(3x+1)\Gamma(x+1)^3}.$$
 Then
 \begin{equation}
 \sum_{k=0}^\infty g_2'(k)=0,\ \sum_{k=0}^\infty g_2''(k)=-48\pi^2
 \ \ \t{and}\ \ \sum_{k=0}^\infty g_2'''(k)=10752\zeta(3).
 \end{equation}
 \end{conjecture}
 \begin{remark} Note that
 $$\sum_{k=0}^\infty g_2(k)=\sum_{k=0}^\infty \f{(74k+7)\bi{2k}k\bi{3k}k\bi{6k}{3k}}{(2k+1)4096^k}=8$$
 by \cite[Example 60]{CZ}.
 \end{remark}
 
 \begin{conjecture} [2026-04-10] For $x>-1/3$, define
 $$g_3(x)=\f{(27x^2+18x+2)\Gamma(3x+1)^2}{(2x+1)729^x\Gamma(x+1)^4\Gamma(2x+1)}.$$
 Then
 \begin{equation}\sum_{k=0}^\infty g_3'(k)=0,\
 \sum_{k=0}^\infty g_3''(k)=-12\sqrt3\,\pi
 \ \ \t{and}\ \ \sum_{k=0}^\infty g_3'''(k)=1456L_{-3}(2).
 \end{equation}
 \end{conjecture}
\begin{remark} Note that
$$\sum_{k=0}^\infty g_3(k)=\sum_{k=0}^\infty\f{(27k^2+18k+2)\bi{2k}k\bi{3k}k^2}{(2k+1)729^k}=\f{9\sqrt3}{2\pi}$$
by \cite[Example 81]{CZ}.
\end{remark}

 \begin{conjecture} [2026-04-10] For $x>-1/4$, define
 $$g_4(x)=\f{(66x^2+37x+4)\Gamma(3x+1)\Gamma(4x+1)}{(2x+1)729^x\Gamma(x+1)^3\Gamma(2x+1)^2}.$$
 Then
 \begin{equation}\sum_{k=0}^\infty g_4'(k)=0,\
 \sum_{k=0}^\infty g_4''(k)=-18\sqrt3\,\pi
 \ \ \t{and}\ \ \sum_{k=0}^\infty g_4'''(k)=2187L_{-3}(2).
 \end{equation}
 \end{conjecture}
\begin{remark} Note that
$$\sum_{k=0}^\infty g_4(k)=\sum_{k=0}^\infty\f{(66k^2+37k+4)\bi{3k}k^2\bi{4k}k}{(2k+1)729^k}=\f{27\sqrt3}{2\pi}$$
as conjectured by the author \cite[(5.21)]{S24} and later confirmed by Au \cite{Au}.
\end{remark}

\begin{conjecture} [2026-04-10] For $x>-1/4$, define
 $$g_5(x)=\f{e^{\pi ix}(258x^2+167x+20)\Gamma(3x+1)\Gamma(4x+1)}{(2x+1)729^x\Gamma(x+1)^3\Gamma(2x+1)^2}.$$
 Then
 \begin{equation}\sum_{k=0}^\infty g_5'(k)=27\sqrt3\, i,\
 \sum_{k=0}^\infty g_5''(k)=-90\sqrt3\,\pi
 \end{equation}
 and \begin{equation} \sum_{k=0}^\infty g_5'''(k)=108\l(81L_{-3}(2)-2\sqrt3\pi^2i\r).
 \end{equation}
 \end{conjecture}
\begin{remark} Note that
$$\sum_{k=0}^\infty g_5(k)=\sum_{k=0}^\infty\f{(258k^2+167k+20)\bi{3k}k^2\bi{4k}k}{(2k+1)(-729)^k}=\f{27\sqrt3}{\pi}$$
by \cite[Example 49]{CZ}.
\end{remark}

\begin{conjecture} [2026-04-10] For $x>-1/4$, define
 $$g_6(x)=\f{(48x^2+32x+3)\Gamma(4x+1)^2}{(2x+1)4096^x\Gamma(x+1)^2\Gamma(2x+1)^3}.$$
 Then
 \begin{equation}\sum_{k=0}^\infty g_6'(k)=0,
 \  \sum_{k=0}^\infty g_6''(k)=-32\sqrt2\,\pi
 \ \ \t{and}\ \ \sum_{k=0}^\infty g_6'''(k)=3072L_{-8}(2).
 \end{equation}
 \end{conjecture}
\begin{remark} Note that
$$\sum_{k=0}^\infty g_6(k)=\sum_{k=0}^\infty\f{(48k^2+32k+3)\bi{2k}k\bi{4k}{2k}^2}{(2k+1)4096^k}=\f{8\sqrt2}{\pi}$$
by \cite[Example 82]{CZ}.
\end{remark}

\begin{conjecture} [2026-04-10] For $x>-1/4$, define
 $$g_7(x)=\f{(296x^2+182x+15)\Gamma(4x+1)\Gamma(6x+1)}{(2x+1)2^{14x}\Gamma(x+1)\Gamma(2x+1)^3\Gamma(3x+1)}.$$
 Then
 \begin{equation}\sum_{k=0}^\infty g_7'(k)=0,
 \  \sum_{k=0}^\infty g_7''(k)=-256\pi
 \ \ \t{and}\ \ \sum_{k=0}^\infty g_7'''(k)=24576G.
 \end{equation}
 \end{conjecture}
\begin{remark} Note that
$$\sum_{k=0}^\infty g_7(k)=\sum_{k=0}^\infty\f{(296k^2+182k+15)\bi{3k}k\bi{4k}{2k}\bi{6k}{3k}}{(2k+1)2^{14k}}=\f{64}{\pi}$$
by \cite[Example 59]{CZ}.
\end{remark}

\begin{conjecture} [2026-04-10] For $x>-1/4$, define
 $$g_8(x)=\f{(88x^3+108x^2+36x+3)\Gamma(4x+1)^2}{(3x+1)(3x+2)1024^x\Gamma(x+1)^3\Gamma(2x+1)\Gamma(3x+1)}.$$
 Then
 \begin{equation}\sum_{k=0}^\infty g_8'(k)=0,
 \  \sum_{k=0}^\infty g_8''(k)=-16\pi
 \ \ \t{and}\ \ \sum_{k=0}^\infty g_8'''(k)=1152G.
 \end{equation}
 \end{conjecture}
\begin{remark} Note that
$$\sum_{k=0}^\infty g_8(k)=\sum_{k=0}^\infty\f{(88k^3+108k^2+36k+3)\bi{3k}k\bi{4k}{k}^2}{(3k+1)(3k+2)1024^k}=\f{8}{\pi}$$
by \cite[Example 110]{CZ}.
\end{remark}

\begin{conjecture} [2026-04-10] For $x>-1/4$, define
 $$g_9(x)=\f{(368x^3+400x^2+118x+9)\Gamma(4x+1)^2}{(3x+1)(3x+2)4096^x\Gamma(x+1)^3\Gamma(2x+1)\Gamma(3x+1)}.$$
 Then
 \begin{equation}\sum_{k=0}^\infty g_9'(k)=0,
 \  \sum_{k=0}^\infty g_9''(k)=-64\pi
 \ \ \t{and}\ \ \sum_{k=0}^\infty g_9'''(k)=4608G.
 \end{equation}
 \end{conjecture}
\begin{remark} Note that
$$\sum_{k=0}^\infty g_9(k)=\sum_{k=0}^\infty\f{(368k^3+400k^2+118k+9)\bi{3k}k\bi{4k}{k}^2}{(3k+1)(3k+2)4096^k}=\f{16}{\pi}$$
by \cite[Example 106]{CZ}.
\end{remark}

\begin{conjecture} [2026-04-10] For $x>-1/4$, define
 $$g_{10}(x)=\f{e^{\pi ix}(896x^3+992x^2+296x+21)\Gamma(4x+1)^2}{(3x+1)(3x+2)2^{14x}\Gamma(x+1)^3\Gamma(2x+1)\Gamma(3x+1)}.$$
 Then
 \begin{equation}\sum_{k=0}^\infty g_{10}'(k)=32i,
 \  \sum_{k=0}^\infty g_{10}''(k)=-256\pi
 \ \ \t{and}\ \ \sum_{k=0}^\infty g_{10}'''(k)=64(288G-11\pi^2i).
 \end{equation}
 \end{conjecture}
\begin{remark} Note that
$$\sum_{k=0}^\infty g_{10}(k)=\sum_{k=0}^\infty\f{(896k^3+992k^2+296k+21)\bi{3k}k\bi{4k}{k}^2}{(3k+1)(3k+2)(-2^{14})^k}=\f{32}{\pi}$$
by \cite[Example 97]{CZ}.
\end{remark}

\begin{conjecture} [2026-04-10] For $x>-1/6$, define
 $$g_{11}(x)=\f{(360x^3+612x^2+230x+15)\Gamma(6x+1)^2}{(3x+1)(3x+2)2^{15x}\Gamma(x+1)\Gamma(2x+1)\Gamma(3x+1)^3}.$$
 Then
 \begin{equation}\sum_{k=0}^\infty g_{11}'(k)=0,
 \  \sum_{k=0}^\infty g_{11}''(k)=-128\sqrt2\,\pi
 \ \ \t{and}\ \ \sum_{k=0}^\infty g_{11}'''(k)=18432L_{-8}(2).
 \end{equation}
 \end{conjecture}
\begin{remark} Note that
$$\sum_{k=0}^\infty g_{11}(k)=\sum_{k=0}^\infty\f{(360k^3+612k^2+230k+15)\bi{3k}k\bi{6k}{3k}^2}{(3k+1)(3k+2)2^{15k}}$$
was conjectured to have the value $32\sqrt2/\pi$ by the author \cite[(5.20)]{II}.
\end{remark}

\begin{conjecture} [2026-04-10] For $x>-1/6$, define
 $$g_{12}(x)=\f{e^{\pi ix}(28x^2+10x+1)\Gamma(2x+1)^6\Gamma(3x+1)}{(6x+1)64^x\Gamma(x+1)^9\Gamma(6x+1)}.$$
 Then
 \begin{gather}\sum_{k=0}^\infty g_{12}'(k)=3i,
 \  \sum_{k=0}^\infty g_{12}''(k)=-12\pi,
 \\ \sum_{k=0}^\infty g_{12}'''(k)=288G-30\pi^2 i,
 \ \sum_{k=0}^\infty \mathrm{Im}(g_{12}^{(4)}(k))=1152\pi G.
 \end{gather}
 \end{conjecture}
\begin{remark} Note that
$$\sum_{k=0}^\infty g_{12}(k)=\sum_{k=0}^\infty\f{(28k^2+10k+1)\bi{2k}k^5}{(6k+1)(-64)^k\bi{6k}{3k}\bi{3k}k}=\f{3}{\pi}$$
as conjectured by  \cite[(5.34)]{II} and recently confirmed by Q.-H. Hou and the author.
\end{remark}

\begin{conjecture} [2026-04-10] For $x>-1/2$, define
 $$g_{13}(x)=\f{e^{\pi i x}(20x^2+8x+1)\Gamma(2x+1)^5}{4096^x\Gamma(x+1)^{10}}.$$
 Then
 \begin{gather}
 \sum_{k=0}^\infty g_{13}''(k)=-16,\ \sum_{k=0}^\infty g_{13}'''(k)=-32\pi i,
 \\\sum_{k=0}^\infty g_{13}^{(4)}(k)=256\pi^2
 \ \ \t{and}\ \ \sum_{k=0}^\infty g_{13}^{(5)}(k)=64\l(17\pi^3i-840\zeta(3)\r).
 \end{gather}
 \end{conjecture}
\begin{remark} Note that
$$\sum_{k=0}^\infty g_{13}(k)=\sum_{k=0}^\infty \f{20k^2+8k+1}{(-4096)^k}\bi{2k}k^5=\f{8}{\pi^2}$$
by Guillera \cite[Indentity 8]{G08}.
Also, we have
$\sum_{k=0}^\infty g_{13}'(k)=8i/\pi$ by the conjectural formulas (141) and (142) of the author \cite{harmonic} confirmed by C. Wei \cite{W}.
\end{remark}

\begin{conjecture} [2026-04-10] For $x>-1/2$, define
 $$g_{14}(x)=\f{e^{\pi i x}(820x^2+180x+13)\Gamma(2x+1)^5}{2^{20x}\Gamma(x+1)^{10}}.$$
 Then
 \begin{gather}
 \sum_{k=0}^\infty g_{14}''(k)=-768,\ \sum_{k=0}^\infty g_{14}'''(k)=-2048\pi i,
 \\\sum_{k=0}^\infty g_{14}^{(4)}(k)=43008\pi^2
 \ \ \t{and}\ \ \sum_{k=0}^\infty g_{14}^{(5)}(k)=1024\l(197\pi^3i-13440\zeta(3)\r).
 \end{gather}
 \end{conjecture}
\begin{remark} Note that
$$\sum_{k=0}^\infty g_{14}(k)=\sum_{k=0}^\infty \f{820k^2+180k+13}{(-2^{20})^k}\bi{2k}k^5=\f{128}{\pi^2}$$
by Guillera \cite[Indentity 9]{G08}.
Also, we have
$\sum_{k=0}^\infty g_{14}'(k)=128i/\pi$ by the conjectural formula (148) of the author \cite{harmonic} confirmed by C. Wei \cite{W}.
\end{remark}

\begin{conjecture} [2026-04-11] For $x>-1/2$, define
$$g_{15}(x)=\f{(77x^2+27x+3)\Gamma(2x+1)^3\Gamma(3x+1)}{2^{12x}\Gamma(x+1)^9}.$$
Then
\begin{gather}\sum_{k=0}^\infty g_{15}'(k)=\sum_{k=0}^\infty g_{15}'''(k)=0,
\ \sum_{k=0}^\infty g_{15}''(k)=-32,\\ \sum_{k=0}^\infty g_{15}^{(4)}(k)=768\pi^2,
\ \sum_{k=0}^\infty g_{15}^{(5)}(k)=-215040\zeta(3).
\end{gather}
\end{conjecture}
\begin{remark} Note that
$$\sum_{k=0}^\infty g_{15}(k)=\sum_{k=0}^\infty\f{74k^2+27k+3}{4096^k}\bi{2k}k^4\bi{3k}{k}=\f{48}{\pi^2}$$
by \cite{G11}.
\end{remark}

\begin{conjecture} [2026-04-10] For $x>-1/2$, define
$$g_{16}(x)=\f{(120x^2+34x+3)\Gamma(2x+1)^2\Gamma(4x+1)}{2^{16x}\Gamma(x+1)^8}.$$
Then
\begin{gather}\sum_{k=0}^\infty g_{16}'(k)=\sum_{k=0}^\infty g_{16}'''(k)=0,
\ \sum_{k=0}^\infty g_{16}''(k)=-64,\\ \sum_{k=0}^\infty g_{16}^{(4)}(k)=2560\pi^2,
\ \sum_{k=0}^\infty g_{16}^{(5)}(k)=-860160\zeta(3).
\end{gather}
\end{conjecture}
\begin{remark} Note that
$$\sum_{k=0}^\infty g_{16}(k)=\sum_{k=0}^\infty\f{120k^2+34k+3}{2^{16k}}\bi{2k}k^4\bi{4k}{2k}=\f{32}{\pi^2}$$
by \cite[Identity 10]{G08}.
\end{remark}

\begin{conjecture} [2026-04-12] For $x>-1/2$, define
$$g_{17}(x)=\f{e^{\pi ix}(252x^2+63x+5)\Gamma(3x+1)\Gamma(4x+1)}{24^{4x}\Gamma(x+1)^7}.$$
Then
\begin{gather}\sum_{k=0}^\infty g_{17}'(k)=0,
\ \sum_{k=0}^\infty g_{17}''(k)=-192,\ \sum_{k=0}^\infty g_{17}'''(k)=-480\pi i,
\\\sum_{k=0}^\infty g_{17}^{(4)}(k)=8448\pi^2,
\ \sum_{k=0}^\infty g_{17}^{(5)}(k)=768\l(51\pi^3i-3640\zeta(3)\r).
\end{gather}
\end{conjecture}
\begin{remark} Note that
$$\sum_{k=0}^\infty g_{17}(k)=\sum_{k=0}^\infty\f{252k^2+63k+5}{(-24^4)^k}\bi{2k}k^2\bi{3k}k^2\bi{4k}{k}=\f{48}{\pi^2}$$
by \cite{CG}.
\end{remark}

\begin{conjecture} [2026-04-12] For $x>-1/6$, define
$$g_{18}(x)=\f{(532x^2+126x+9)\Gamma(6x+1)}{10^{6x}\Gamma(x+1)^6}.$$
Then
\begin{gather}\sum_{k=0}^\infty g_{18}'(k)=\sum_{k=0}^\infty g_{18}'''(k)=0,
\ \sum_{k=0}^\infty g_{18}''(k)=-250,\\ \sum_{k=0}^\infty g_{18}^{(4)}(k)=15000\pi^2,
\ \sum_{k=0}^\infty g_{18}^{(5)}(k)=-6510000\zeta(3).
\end{gather}
\end{conjecture}
\begin{remark} Note that
$$\sum_{k=0}^\infty g_{18}(k)=\sum_{k=0}^\infty\f{532k^2+126k+9}{10^{6k}}\bi{2k}k^2\bi{3k}k^2\bi{6k}{3k}=\f{375}{4\pi^2}$$
in view of \cite{CG}.
\end{remark}

\begin{conjecture} [2026-04-12] For $x>-1/6$, define
$$g_{19}(x)=\f{e^{\pi ix}(1930x^2+549x+45)\Gamma(6x+1)}{2^{18x}\Gamma(x+1)^6}.$$
Then
\begin{gather}\sum_{k=0}^\infty g_{19}'(k)=\f{384i}{\pi},
\ \sum_{k=0}^\infty g_{19}''(k)=-1024,\ \sum_{k=0}^\infty g_{19}'''(k)=-2304\pi i,
\\ \sum_{k=0}^\infty g_{19}^{(4)}(k)=36864\pi^2,\
\t{and}\ \sum_{k=0}^\infty g_{19}^{(5)}(k)=2048\l(83\pi^3 i-6720\zeta(3)\r).
\end{gather}
\end{conjecture}
\begin{remark} Note that
$$\sum_{k=0}^\infty g_{19}(k)=\sum_{k=0}^\infty\f{1930k^2+549k+45}{(-2^{18})^k}\bi{2k}k^2\bi{3k}k^2\bi{6k}{3k}$$
has the value $384/\pi^2$ as proved by Au \cite{Au}.
\end{remark}

\begin{conjecture} [2026-04-12] For $x>-1/6$, define
$$g_{20}(x)=\f{e^{\pi ix}(5418x^2+693x+29)\Gamma(6x+1)}{(2^{18}3^65^3)^k\Gamma(x+1)^6}.$$
Then
\begin{gather}\sum_{k=0}^\infty g_{20}'(k)=\f{128\sqrt5\,i}{\pi},
\ \sum_{k=0}^\infty g_{20}''(k)=-2048\sqrt5,\\ \sum_{k=0}^\infty g_{20}'''(k)=-5888\sqrt5\pi i,
\ \sum_{k=0}^\infty g_{20}^{(4)}(k)=352256\sqrt5\pi^2.
\end{gather}
\end{conjecture}
\begin{remark} Note that
$$\sum_{k=0}^\infty g_{20}(k)=\sum_{k=0}^\infty\f{5418k^2+693k+29}{(-2^{18}3^65^3)^k}\bi{2k}k^2\bi{3k}k^2\bi{6k}{3k}$$
was conjectured to have the value $128\sqrt5/\pi^2$ (cf. \cite{CG}).
\end{remark}

\begin{conjecture} [2026-04-12] For $x>-1/6$, define
$$g_{21}(x)=\f{(1640x^2+278x+15)e^{\pi ix}\Gamma(4x+1)\Gamma(6x+1)}{(2^{22}3^3)^x\Gamma(x+1)^5\Gamma(2x+1)\Gamma(3x+1)}.$$
Then
\begin{gather}\sum_{k=0}^\infty g_{21}'(k)=\f{256i}{\sqrt3\,\pi},
\ \sum_{k=0}^\infty g_{21}''(k)=-\f{2048}{\sqrt3},
\\\sum_{k=0}^\infty g_{21}''(k)=-\f{5632\pi i}{\sqrt3}
\ \t{and}\ \sum_{k=0}^\infty g_{21}^{(4)}(k)=\f{204800}{\sqrt3}\pi^2.
\end{gather}
\end{conjecture}
\begin{remark} Note that
$$\sum_{k=0}^\infty g_{21}(k)=\sum_{k=0}^\infty
\f{1640k^2+278k+15}{(-2^{22}3^3)^k}\bi{2k}k^2\bi{3k}k\bi{4k}{2k}\bi{6k}{3k}$$
was conjectured to have the value $256/(\sqrt3\pi^2)$ (cf. \cite{CG}).
\end{remark}

\begin{conjecture} [2026-04-12] For $x>-1/8$, define
$$g_{22}(x)=(1920x^2+304x+15)\f{\Gamma(2x+1)\Gamma(8x+1)}{\Gamma(x+1)^6\Gamma(4x+1)}.$$
Then
\begin{gather}\sum_{k=0}^\infty g_{22}'(k)=\sum_{k=0}^\infty g_{22}'''(k)=0,
\\\sum_{k=0}^\infty g_{22}''(k)=-448\sqrt{7}
\ \t{and}\ \sum_{k=0}^\infty g_{22}^{(4)}(k)=55552\sqrt7\pi^2.
\end{gather}
\end{conjecture}
\begin{remark} Note that
$$\sum_{k=0}^\infty g_{22}(k)=\sum_{k=0}^\infty
\f{1920k^2+304k+15}{(2^{18}7^4)^k}\bi{2k}k^3\bi{4k}{2k}\bi{8k}{4k}$$
was conjectured to have the value $56\sqrt7/\pi^2$ (cf. \cite{CG}).
\end{remark}

\begin{conjecture} [2026-04-13] For $x>-1/4$, define
$$g_{23}(x)=\f{(1344x^3+944x^2+156x+9)\Gamma(4x+1)^3}{(2x+1)2^{22x}\Gamma(x+1)^{4}\Gamma(2x+1)^4}.$$
Then
\begin{gather}\sum_{k=0}^\infty g_{23}'(k)=\sum_{k=0}^\infty g_{23}'''(k)=0,
\\ \sum_{k=0}^\infty g_{23}''(k)=-256\sqrt2,\ \sum_{k=0}^\infty g_{23}^{(4)}(k)=19456\sqrt2\pi^2.
\end{gather}
\end{conjecture}
\begin{remark} Note that
$$\sum_{k=0}^\infty g_{23}(k)=\sum_{k=0}^\infty\f{(1344k^3+944k^2+156k+9)\bi{2k}k^2\bi{4k}{2k}^3}{(2k+1)2^{22k}}$$
was conjectured to have the value $64\sqrt2/\pi^2$ by the author \cite[(5.31)]{II}.
\end{remark}

\begin{conjecture} [2026-04-13] For $x>-1/6$, define
$$g_{24}(x)=\f{(92x^3+54x^2+12x+1)\Gamma(2x+1)^8\Gamma(3x+1)}{(6x+1)256^x\Gamma(x+1)^{13}\Gamma(6x+1)}.$$
Then
\begin{gather}\sum_{k=0}^\infty g_{24}'(k)=\sum_{k=0}^\infty g_{24}'''(k)=0,
\ \sum_{k=0}^\infty g_{24}''(k)=-12,\\ \sum_{k=0}^\infty g_{24}^{(4)}(k)=192\pi^2,
\ \sum_{k=0}^\infty g_{24}^{(5)}(k)=-30240\zeta(3).
\end{gather}
\end{conjecture}
\begin{remark} Note that
$$\sum_{k=0}^\infty g_{24}(k)=\sum_{k=0}^\infty\f{(92k^3+54k^2+12k+1)\bi{2k}k^7}{(6k+1)256^k\bi{6k}{3k}\bi{3k}k}=\f{12}{\pi^2}$$
as conjectured by the author \cite[(5.30)]{S24} and confirmed by Au \cite{Au}.
\end{remark}

\begin{conjecture} [2026-04-12] For $x>-1/2$, define
$$g_{25}(x)=(168x^3+76x^2+14x+1)\f{\Gamma(2x+1)^7}{2^{20x}\Gamma(x+1)^{14}}.$$
Then
\begin{gather}\sum_{k=0}^\infty g_{25}'(k)=\sum_{k=0}^\infty g_{25}'''(k)=\sum_{k=0}^\infty g_{25}^{(5)}(k)=0,
\\\sum_{k=0}^\infty g_{25}''(k)=-\f{64}{\pi},\ \sum_{k=0}^\infty g_{25}^{(4)}(k)=1024\pi,
\ \sum_{k=0}^\infty g_{25}^{(5)}(k)=-131584\pi^3.
\end{gather}
\end{conjecture}
\begin{remark} Note that
$$\sum_{k=0}^\infty g_{25}(k)=\sum_{k=0}^\infty(168k^3+76k^2+14k+1)\f{\bi{2k}k^7}{2^{20k}}=\f{32}{\pi^3}$$
as conjectured by B. Gourevich and confirmed by Au \cite{Au}.
\end{remark}

\begin{conjecture} [2026-04-12] For $x>-1/2$, define
$$g_{26}(x)=\f{e^{\pi ix}P(x)\Gamma(2x+1)^4\Gamma(3x+1)\Gamma(4x+1)}{2^{24x}\Gamma(x+1)^{15}},$$
where $P(x)=4528x^4+3180x^3+972x^2+147x+9$.
Then
\begin{gather}\sum_{k=0}^\infty g_{26}'(k)=\f{768i}{\pi^3},\
\sum_{k=0}^\infty g_{26}''(k)=-\f{1536}{\pi^2},
\ \sum_{k=0}^\infty g_{26}'''(k)=-\f{3072i}{\pi},
\\\sum_{k=0}^\infty g_{26}^{(4)}(k)=12288,\ \sum_{k=0}^\infty g_{25}^{(5)}(k)=43008\pi i,
\\ \sum_{k=0}^\infty g_{26}^{(6)}(k)=-454656\pi^2,\ \sum_{k=0}^\infty g_{25}^{(7)}(k)=-2629632\pi^3i,
\\ \sum_{k=0}^\infty g_{26}^{(8)}(k)=130613248\pi^4,\ \sum_{k=0}^\infty g_{25}^{(9)}(k)=
196608\l(5671\pi^5i-5624640\zeta(5)\r).
\end{gather}
\end{conjecture}
\begin{remark} Note that
$$\sum_{k=0}^\infty g_{26}(k)=\sum_{k=0}^\infty\f{P(k)}{(-2^{24})^k}\bi{2k}k^7\bi{3k}k\bi{4k}{2k}$$
was conjectured to have the value $768/\pi^4$ (cf. \cite{CG}).
\end{remark}

\begin{conjecture} [2026-04-12] For $x>-1/2$, define
$$g_{27}(x)=\f{Q(x)\Gamma(2x+1)^6\Gamma(4x+1)}{2^{32x}\Gamma(x+1)^{16}},$$
where $Q(x)=43680x^4+20632x^3+4340x^2+466x+21$.
Then
\begin{gather}\sum_{k=0}^\infty g_{27}'(k)=\sum_{k=0}^\infty g_{27}'''(k)=\sum_{k=0}^\infty g_{27}^{(5)}(k)
=\sum_{k=0}^\infty g_{27}^{(7)}(k)=0,
\\
\sum_{k=0}^\infty g_{27}''(k)=-\f{8192}{\pi^2},
\ \sum_{k=0}^\infty g_{27}^{(4)}(k)=180224,
\\\sum_{k=0}^\infty g_{27}^{(6)}(k)=-14876672\pi^2,\ \sum_{k=0}^\infty g_{27}^{(8)}(k)=7068975104\pi^4,
\\\sum_{k=0}^\infty g_{27}^{(9)}(k)=-70774350151680\zeta(5).
\end{gather}
\end{conjecture}
\begin{remark} Note that
$$\sum_{k=0}^\infty g_{27}(k)=\sum_{k=0}^\infty\f{Q(k)}{2^{32k}}\bi{2k}k^8\bi{4k}{2k}$$
was conjectured to have the value $2048/\pi^4$ (cf. \cite{CG}).
\end{remark}

\end{document}